\documentclass[preprint,12pt]{elsarticle}

\usepackage{amssymb}
\newcommand\sbullet[1][.5]{\mathbin{\vcenter{\hbox{\scalebox{#1}{$\bullet$}}}}}
\usepackage{amsthm}
\usepackage{amsmath}
\usepackage{graphics}
\usepackage{epsfig}
\usepackage{amssymb}
\usepackage{hyperref}
\usepackage{graphicx}
\usepackage{subfigure}
\usepackage{color}
\usepackage{tikz}
\usetikzlibrary{patterns}
\usepackage[ruled,vlined]{algorithm2e}
\usepackage{hyperref}
\usepackage{a4wide}
\usepackage{amsmath}
\usepackage{amsfonts}
\usepackage{enumerate}
\newcommand{\pf}{\noindent {\bf Proof: }}

\newtheorem*{theoremaux}{Theorem \theoremauxnum}
\gdef\theoremauxnum{1}

\newtheorem{lemma}{\bf Lemma}[section]

\newtheorem{theorem}{\bf Theorem}[section]

\newtheorem{corollary}[lemma]{\bf Corollary}
\newtheorem{definition}{\bf Definition}[section]

\journal{~}

\begin{document}
	
	\begin{frontmatter}
		
		
			
\author[1]{Angsuman Das\corref{cor1}}
\ead{angsuman.maths@presiuniv.ac.in}

\author[1]{Arnab Mandal}
\ead{arnab.maths@presiuniv.ac.in}

\address[1]{Department of Mathematics,\\ Presidency University, Kolkata, India} 
\cortext[cor1]{Corresponding author}
		
\title{Solvability of a group based on its number of subgroups}

\begin{abstract}
In this paper, we provide some conditions of (super)-solvability and nilpotency of a finite group $G$ based on its number of subgroups $Sub(G)$.	Our results generalize the classification of finite groups with less than $20$ subgroups by Betz and Nash. We also provide an application of our results in studying comaximal subgroup graph of a group. Finally, we conclude with some open issues.
\end{abstract}
		
\end{frontmatter}

\section{Introduction} 
A major research area in finite group theory is to identify a group $G$ from partial information about it. For instances, the order profile/sequence of its elements \cite{hiranya-cameron}, sum of orders of its elements \cite{amiri}, \cite{asad}, \cite{herzog}, the number of subgroups \cite{betz-nash},  graphs defined on groups \cite{comaximal-1}, \cite{comaximal-2} etc. have been used to predict the nature of the underlying group. 

In this paper, we focus on the number of subgroups, denoted by $Sub(G)$ of a finite group $G$, which includes the trivial subgroup and the group $G$ itself. Classifying groups by the number of subgroups is a classically interesting problem and dates back to 1939 when Miller in a series of papers \cite{miller1,miller2,miller3,miller4,miller5} classified groups $G$ with $Sub(G)\leq 16$. Slattery \cite{slattery} in 2016 studied the same problem with a different approach of {\it similarity of groups}. Recently, Betz and Nash (2022) \cite{betz-nash} corrected and extended these results by classifying all abelian groups $G$ with $Sub(G)\leq 22$ and all non-abelian groups $G$ with $Sub(G)\leq 19$. In fact, they \cite{betz-nash-2} also extended their classification for abelian groups upto $Sub(G)\leq 49$. The list of groups $G$ with $Sub(G)\leq 11$ is given below (Table \ref{list<12}). For the full table, consisting of groups $G$ with $Sub(G)\leq 19$, the readers are referred to Table 3 in \cite{betz-nash}. For the number of groups with $Sub(G)=k$, where $1\leq k \leq 19$, one can refer to OEIS A274847 \cite{oeis}.

\begin{table}
	\centering
	\begin{center}
		\begin{tabular}{|c|c|}
			\hline $Sub(G)$ & Groups ($G$) \\
			\hline \hline \hline  $1$ & Trivial group \\
			\hline  $2$ & $\mathbb{Z}_{p}$ \\
			\hline  $3$ & $\mathbb{Z}_{p^2}$ \\
			\hline $4$ & $\mathbb{Z}_{p^3}$ or $\mathbb{Z}_{pq}$\\
			\hline $5$ & $\mathbb{Z}_{p^4}$ or  $\mathbb{Z}_2\times \mathbb{Z}_2$\\
			\hline $6$ & $\mathbb{Z}_{p^5}$, $\mathbb{Z}_{p^2q}$, $\mathbb{Z}_2\times \mathbb{Z}_2$, $S_3$, $Q_8$.\\
			\hline $7$ & $\mathbb{Z}_{p^6}$\\
			\hline $8$ & $\mathbb{Z}_{p^7}$, $\mathbb{Z}_{p^3q}$, $\mathbb{Z}_{pqr}$, $\mathbb{Z}_4\times \mathbb{Z}_2$, $\mathbb{Z}_5\times \mathbb{Z}_5$, $Dic_{12}$, $D_5$.\\
			\hline $9$ & $\mathbb{Z}_{p^8}$, $\mathbb{Z}_{p^2q^2}$\\
			\hline $10$ & $\mathbb{Z}_{p^9}$, $\mathbb{Z}_{p^4q}$, $\mathbb{Z}_2\times \mathbb{Z}_2\times \mathbb{Z}_p$, $\mathbb{Z}_7 \times \mathbb{Z}_7$, $\mathbb{Z}_9 \times \mathbb{Z}_3$,\\  
			& $\mathbb{Z}_7\rtimes \mathbb{Z}_3$, $\mathbb{Z}_3\rtimes \mathbb{Z}_8$, $D_4$, $D_7$, $Dic_{20}$, $A_4$\\
			\hline $11$ & $\mathbb{Z}_{p^{10}}$, $\mathbb{Z}_8\times \mathbb{Z}_2$, $Q_{16}$, $M_{16}$\\
			\hline
		\end{tabular}
		\caption{Classification of Groups by $Sub(G)$}
		\label{list<12}  
	\end{center}
	
\end{table}

As mentioned by the authors in \cite{betz-nash}, since their technique uses GAP, it is computationally difficult to extend their results beyond $19$ subgroups, especially for the non-abelian case. Keeping this in mind, instead of determining the exact group, we focus on finding the nature of the groups $G$, e.g., solvability, supersolvability, nilpotency etc, from a given value of $Sub(G)$. It is to be noted that we also use GAP in our methods, but due to some tight bounding arguments we manage to reduce the search space for GAP.

\subsection{Preliminaries and Motivation}
As evident from Table \ref{list<12} (above) and Table 1 \& 3 in \cite{betz-nash}, we observe that 
\begin{itemize}
	\item if $Sub(G)=1,2,3,4,7,9,13$, then $G$ is cyclic.
	\item if $Sub(G)\leq 5$ or $Sub(G)=7,9,13$, then $G$ is abelian.
	\item if $Sub(G)\leq 5$ or $Sub(G)=7,9,11,13,17$, then $G$ is nilpotent.
	\item if $Sub(G)\leq 19$ with $Sub(G)\neq 10,15$, then $G$ is supersolvable.
	\item if $Sub(G)\leq 19$, then $G$ is solvable.
\end{itemize}
Moreover, the only non-supersolvable groups $G$ with $Sub(G)\leq 19$ are $A_4$, $SL(2,3)$ and $(\mathbb{Z}_2\times \mathbb{Z}_2)\rtimes \mathbb{Z}_9$ with $Sub(A_4)=10$  and $Sub(SL(2,3))=Sub((\mathbb{Z}_2\times \mathbb{Z}_2)\rtimes \mathbb{Z}_9)=15$. This motivates us to define the following:  
\begin{definition}
	A positive integer $n$ is called $\sbullet[1.25]$ enforcing number, where $\sbullet[1.25] \in \{ $Solvable, Supersolvable, Nilpotent, Abelian, Cyclic$\}$, if $Sub(G)=n$ implies that $G$ is a $\sbullet[1.25]$ group. The set of all $\sbullet[1.25]$ enforcing numbers is called the $\sbullet[1.25]$ enforcing set and is denoted by $\mathsf{Enf}(S),\mathsf{Enf}(SS)$, $\mathsf{Enf}(N)$, $\mathsf{Enf}(A)$ and $\mathsf{Enf}(C)$ respectively.
\end{definition}
It is clear from the definition that $\mathsf{Enf}(C)\subseteq \mathsf{Enf}(A) \subseteq \mathsf{Enf}(N)\subseteq \mathsf{Enf}(SS)\subseteq \mathsf{Enf}(S)$. Note that all the inclusions are proper and the numbers $1,2,\ldots,19$ are fully classified with respect to $\sbullet[1.25]$ enforcing numbers. A solvable enforcing number is called a {\it strictly solvable enforcing number} if it is not a supersolvable enforcing number. We also observe that if $n$ is not a $\sbullet[1.25]$ enforcing number, then any multiple of $n$ is also not a $\sbullet[1.25]$ enforcing number, because if $G$ is a non-$\sbullet[1.25]$ group with $Sub(G)=n$ and $p$ is a prime such that $p\nmid |G|$, then $H=G \times \mathbb{Z}_{p^{k-1}}$ is also a non-$\sbullet[1.25]$ group with $Sub(H)=nk$. Hence $20$ is not a supersolvable enforcing number. Again as $Sub(Q_8\rtimes \mathbb{Z}_9)=21$ and $Q_8\rtimes \mathbb{Z}_9$ is not supersolvable, $21$ is also not a supersolvable enforcing number. Thus $10,15,20,21 \not\in \mathsf{Enf}(SS)$.

\subsection{Our Contribution}
In this article, we show that \begin{enumerate}
	\item If $1\leq Sub(G)\leq 76$ and $Sub(G)\neq 59,76$, then $G$ is solvable. (Section \ref{solv-enf-numbers}.)
 \item If $G$ is a non-solvable group with 
 \begin{itemize}
     \item $Sub(G)=59$, then $G\cong A_5$. (Theorem \ref{Sub=59-A_5})
     \item $Sub(G)=76$, then $G\cong SL(2,5)$. (Theorem \ref{Sub=76-SL(2,5)})
 \end{itemize}
	\item If $22\leq Sub(G)\leq 24$, then $G$ is supersolvable. (Section \ref{ss-enf-numbers})
	\item If $Sub(G)=23$, then $G$ is a $p$-group and hence nilpotent. 
 In fact, we show that if $Sub(G)=23$, then $G$ is one of the $7$ $p$-groups given in Theorem \ref{Sub(G)=23-main-theorem}.
\end{enumerate}
In light of the above results, the previously known results by Betz \& Nash together with GAP computations, we get the first few elements of the enforcing sets, which we demonstrate in the Figure \ref{fig:enforcing_numbers}. Here numbers in \textcolor{blue}{\bf blue} font are cyclic enforcing numbers, \textcolor{cyan}{\bf sky blue} font denotes strictly nilpotent enforcing numbers, \textcolor{orange}{\bf orange} font denotes strictly supersolvable enforcing numbers, {\bf black} font denotes strictly solvable enforcing numbers and numbers with in \boxed{\bf box} are the speculated ones (conjectures), as suggested by GAP after exhaustive search on all finite groups of order $\leq 255$. And the \textcolor{red}{\bf red} font denotes the numbers which are not solvable enforcing, i.e., there exists finite non-solvable groups with those many subgroups.
\begin{figure}
    \centering
$$\begin{array}{cccccccccc}
   \textcolor{blue}{\mathbf{1}}  & \textcolor{blue}{\mathbf{2}}  & \textcolor{blue}{\mathbf{3}} & \textcolor{blue}{\mathbf{4}} & \textcolor{cyan}{\mathbf{5}} & \textcolor{orange}{\mathbf{6}} & \textcolor{blue}{\mathbf{7}} & \textcolor{orange}{\mathbf{8}} & \textcolor{blue}{\mathbf{9}} & \mathbf{10} \\
   \textcolor{cyan}{\mathbf{11}} & \textcolor{orange}{\mathbf{12}} & \textcolor{blue}{\mathbf{13}} & \textcolor{orange}{\mathbf{14}} & \mathbf{15} & \textcolor{orange}{\mathbf{16}} & \textcolor{cyan}{\mathbf{17}} & \textcolor{orange}{\mathbf{18}} & \textcolor{orange}{\mathbf{19}} & \mathbf{20} \\ \mathbf{21} & \textcolor{orange}{\mathbf{22}} & \textcolor{cyan}{\mathbf{23}} & \textcolor{orange}{\mathbf{24}} & \mathbf{25} & \mathbf{26} & \mathbf{27} & \boxed{\textcolor{orange}{\mathbf{28}}} & \boxed{\textcolor{cyan}{\mathbf{29}}} & \mathbf{30} \\ \boxed{\textcolor{orange}{\mathbf{31}}} & \boxed{\textcolor{orange}{\mathbf{32}}} & \boxed{\textcolor{cyan}{\mathbf{33}}} & \mathbf{34} & \mathbf{35} & \mathbf{36} & \mathbf{37} & \mathbf{38} & \boxed{\textcolor{cyan}{\mathbf{39}}} & \mathbf{40} \\
   \mathbf{41} & \mathbf{42} & \boxed{\textcolor{orange}{\mathbf{43}}} & \mathbf{44} & \mathbf{45} & \boxed{\textcolor{orange}{\mathbf{46}}} & \boxed{\textcolor{cyan}{\mathbf{47}}} & \mathbf{48} & \boxed{\textcolor{cyan}{\mathbf{49}}} & \mathbf{50} \\
     \mathbf{51} & \mathbf{52} & \mathbf{53} & \mathbf{54} & \mathbf{55} & \mathbf{56} & \boxed{\textcolor{orange}{\mathbf{57}}} & \mathbf{58} & \textcolor{red}{\mathbf{59}} & \mathbf{60} \\ 
     \boxed{\textcolor{cyan}{\mathbf{61}}} & \boxed{\textcolor{orange}{\mathbf{62}}} & \mathbf{63} & \boxed{\textcolor{orange}{\mathbf{64}}} & \boxed{\textcolor{orange}{\mathbf{65}}} & \mathbf{66} & \boxed{\textcolor{cyan}{\mathbf{67}}} & \mathbf{68} & \mathbf{69} & \mathbf{70} \\ \boxed{\textcolor{orange}{\mathbf{71}}} & \mathbf{72} & \mathbf{73} & \mathbf{74} & \mathbf{75} & \textcolor{red}{\mathbf{76}} 
\end{array}$$    
    \caption{Enforcing Numbers $\leq 76$}
    \label{fig:enforcing_numbers}
\end{figure}

As an application, we demonstrate how our results can be used to study groups from their comaximal subgroup graphs \cite{comaximal-2}. Finally, we conclude with some open questions.

Before going into the main results, we recall some standard results from finite group theory and prove some interesting propositions, which will be repeatedly used in the forthcoming sections.

\begin{itemize}
	\item If $G$ and $H$ are two finite groups with $gcd(|G|,|H|)=1$, then $Sub(G\times H)=Sub(G)\cdot Sub(H)$.
        \item If $N\lhd G$, then $Sub(G)\geq Sub(N)+Sub(G/N)-1$.
        \item If $Sub(G)$ is prime and $Sub(G)\in \mathsf{Enf}(N)$, then $G$ is a $p$-group. (It follows from the fact that nilpotent groups are direct product of its Sylow subgroups.)
        \item (Theorem A \& B, \cite{min-no-of-sg}) If $G$ is a non-cyclic group of order $p^k$ with $k\geq 2$, then 
	$$Sub(G)\geq \left\lbrace \begin{array}{ll}
	6, & \mbox{ if }k=3 \mbox{ and }p=2 \\
	(k-1)p+(k+1), & \mbox{ otherwise.}
	\end{array} \right.$$ 
\end{itemize}

\section{Solvable Enforcing Numbers}\label{solv-enf-numbers}
In this section, we prove that if $1\leq Sub(G)\leq 76$ and $Sub(G)\neq 59,76$, then $G$ is solvable. We also classify the unique non-solvable groups with exactly $59$ or $76$ subgroups.
\begin{theorem}\label{Sub<59-solvable}
Let $G$ be a group such that $Sub(G)<59$, then $G$ is solvable.
\end{theorem}
\pf If possible, let $G$ be a minimum counter-example, i.e., $G$ is a non-solvable group of minimum order such that $Sub(G)<59$. We first show that $G$ must be simple.

We start by noting that all proper subgroups of $G$ are solvable, as if $H$ is a proper subgroup of $G$ which is not solvable, then $H$ and hence $G$ has at least $59$ subgroups, a contradiction. If $H$ is a proper normal subgroup of $G$ and since $H$ is solvable, then $G/H$ must be non-solvable. Now by minimality of $|G|$, we have $Sub(G/H)\geq 59$ and hence by correspondence theorem of subgroups, $Sub(G)\geq 59$, a contradiction. Thus $G$ has no proper normal subgroup, proving that $G$ is simple.

Among all the finite simple groups, $A_5$ has the least number of subgroups and $Sub(A_5)=59$. Thus the theorem holds.\qed 

\begin{theorem}\label{Sub=59-A_5}
Let $G$ be a non-solvable group such that $Sub(G)=59$, then $G\cong A_5$.
\end{theorem}
\pf Since $G$ is non-solvable, there exists a minimal simple group arising as a subquotient of $G$, i.e., $G$ has subgroups $H,N$ such that $N\lhd H$
and $H/N$ is a minimal simple group. If $Sub(H/N)\geq 60$, then $Sub(G)\geq 60$, a contradiction. Thus $Sub(H/N)=59$. If $H$ is a proper subgroup of $G$, then $Sub(G)\geq Sub(H)+1\geq Sub(H/N)+1=60$, a contradiction. Thus $G=H$, i.e., $Sub(G/N)=59$. If $N$ is a non-trivial subgroup of $G$, then $Sub(G)\geq Sub(G/N)+1=60$, again a contradiction. Thus $N$ is the trivial subgroup, and hence $G/N\cong G$ is a minimal simple group. Now from classification of finite minimal simple groups, we get $G\cong A_5$.\qed

\begin{theorem}\label{60<SubG<75}
Let $G$ be a group such that $60\leq Sub(G)\leq 75$, then $G$ is solvable.
\end{theorem}
\pf Let, if possible, $G$ be a non-solvable group of minimum order with $60\leq Sub(G)\leq 75$. We claim that $G$ is simple.

If not, let $H$ be a non-trivial proper normal subgroup of $G$. Then either $H$ or $G/H$ is non-solvable.

{\bf Case 1:} ($H$ is non-solvable.) By minimality of $|G|$, we get $Sub(H)<60$ or $Sub(H)>75$. As $Sub(G)>Sub(H)$, the second case can not occur. So, by Theorem \ref{Sub<59-solvable} and \ref{Sub=59-A_5}, we have $Sub(H)=59$ and $H\cong A_5$.

Also, note that $G/H\cong G/A_5$ is non-trivial group. Hence $G/A_5$ has a subgroup $K/A_5$ of prime order, say $p$, where $K$ is a non-solvable subgroup of $G$ containing $A_5$ and $|K|=60p$. If $p=2,3$ or $5$, it can be checked using GAP that all non-solvable groups of order $120,180$ and $300$ has at least $76$ subgroups. So $p>5$.

Let $L$ be a subgroup of order $p$ in $K$. Then the number of Sylow $p$-subgroups of $K$ is $n_p=1+pk|60$. If $n_p\neq 1$, the possible choices of $p$ are $7,11$ and $19$. Again, it can be exhaustively checked that such case can not occur. So $n_p=1$ and $L\lhd K$. Also, we have $H\cong A_5\lhd K$. Thus $K\cong A_5\times L$ and since $gcd(|A_5|,|L|)=1$, we have $Sub(K)=Sub(A_5)\cdot Sub(L)=118$, which exceeds $Sub(G)$, a contradiction.

{\bf Case 2:} ($G/H$ is non-solvable.) Using similar argument as above, we get $Sub(G/H)=59$, $G/H\cong A_5$ and $|G|=60|H|$. If $|H|=2,3,4$ or $5$, then $|G|=120,180,240$ or $300$. It can be checked using GAP that all non-solvable groups of these orders have at least $76$ subgroups. Thus $|H|\neq 2,3,4,5$.

Let $m \in \{2,3,4,5\}$. Then, by Sylow's theorem, $G$ has a subgroup $K$ of order $m$. If $K$ is a unique subgroup of order $m$ in $G$, then $K\lhd G$, $G/K$ is non-solvable and $G/K\cong A_5$. But this imply that $|H|=|K|=m$, a contradiction. Thus $K$ is not a unique subgroup of order $m$ in $G$ and $H \not\subset K$.

Hence, $G$ at least $(1+2)+(1+2)+(1+3)+(1+5)=16$ non-trivial subgroups (Theorem 5.4.10, \cite{gorenstein}) and none of them  contains $H$. Thus $Sub(G)\geq 1+16+Sub(G/H)=76$, a contradiction.

Combining both the cases, we see that $G$ has no non-trivial proper normal subgroup, i.e., $G$ is simple. As there is no simple group $G$ with $60\leq Sub(G)\leq 75$, the theorem follows.\qed

\begin{theorem}\label{Sub=76-SL(2,5)}
Let $G$ be a non-solvable group such that $Sub(G)=76$, then $G\cong SL(2,5)$.
\end{theorem}
\pf From the classification of finite simple groups, it follows that $G$ is not simple. \\
{\bf Claim 1:} $G$ is a perfect group.\\
{\it Proof of Claim 1:} If not, let $\{e\}\subsetneq G'\subsetneq G$. Since $G/G'$ is abelian, it follows that $G'$ is non-solvable. As $Sub(G')<Sub(G)$, from Theorems \ref{Sub<59-solvable}, \ref{Sub=59-A_5} and \ref{60<SubG<75}, we have $Sub(G')=59$, $G'\cong A_5$ and $G/G'\cong G/A_5$ is a non-trivial group. Let $K/A_5$ be a subgroup of prime order $p$ of $G/A_5$. If $p=2,3$ or $5$, then $K$ is a non-solvable group of order $120,180$ or $300$. Among the non-solvable groups of these orders, it can be checked using GAP, that only tenable option, i.e., having $Sub(G)=76$ is $K\cong SL(2,5)$ of order $120$. However as $SL(2,5)$ has no subgroup isomorphic to $A_5$, we get a contradiction. Thus $p>5$.

Thus $|K|=60p$. Let $n_p$ denote the number of Sylow $p$-subgroups of $K$. If $n_p\neq 1$, then $1<n_p=1+pl|60$. This is possible if $p=7,11$ or $19$. Again, among the non-solvable groups of order $420,660$ and $1140$, it can be checked using GAP, that none of them satisfies $Sub(G)=76$. Thus $n_p=1$, i.e., if $L$ is a Sylow $p$-subgroup of $K$, then $L\lhd K$. Also we have $A_5\lhd K$. Thus $K\cong A_5\times L\cong A_5\times \mathbb{Z}_p$ and $Sub(K)=59\cdot 2=118>76$, a contradiction. Thus Claim 1 holds.

Since $G$ is not simple, $G$ has a non-trivial proper subgroup $N$.\\
{\bf Claim 2:} $N$ is a solvable group.\\
{\it Proof of Claim 2:} If $N$ is not solvable, then by Theorems \ref{Sub<59-solvable}, \ref{Sub=59-A_5} and \ref{60<SubG<75}, we have $Sub(N)=59$ and $N\cong A_5\lhd G$. Let $K/N$ be a subgroup of prime order $p$ in $G/N$. Then $|K|=60p$. Now arguing as in the proof of Claim 1, we get a contradiction. Thus Claim 2 holds.

Since $N$ is solvable, it follows that $G/N$ is non-solvable and hence $G/N\cong A_5$. If there exists a proper subgroup $M$ of $N$, which is normal in $G$, then by similar argument, we get $G/M\cong A_5$. Thus $M=N$. Hence $N$ is a minimal normal subgroup of $G$ and if $M\lhd G$, then $|M|=|N|$.\\
{\bf Claim 3:} $N$ is the unique non-trivial proper normal subgroup of $G$.\\
{\it Proof of Claim 3:} Suppose $M$ is another non-trivial proper normal subgroup of $G$. Then $|M|=|N|$. Moreover $MN$ is a normal subgroup of $G$ with $|MN|>|M|$. Thus $MN=G$, i.e., $$60|N|=|G|=|MN|=\dfrac{|M||N|}{|M\cap N|}=\dfrac{|N|^2}{|M\cap N|}, \mbox{ i.e., }|M\cap N|=|N|/60.$$

Since $N$ is a minimal subgroup, we have $|M\cap N|=1$, i.e., $|N|=60$ and $|G|=3600$. Now, using GAP, one can check that among all perfect groups $G$ of order $3600$, none of them satisfies $Sub(G)=76$. ({\it Note that if $G$ is not known to be perfect, the number of non-solvable groups of order $3600$ is to big to handle with GAP for an exhaustive search}) Thus Claim 3 holds.

 One can also observe that $N$ is characteristically simple, because if not, let $\{e\}\neq M\neq N$ be a characteristic subgroup of $N$. Then $M\lhd G$, a contradiction.

 Now, as $N$ is solvable and characteristically simple, $N$ is the direct product of isomorphic abelian simple groups, i.e., $N\cong \mathbb{Z}^n_p$.

 If $n\geq 4$, then $Sub(N)\geq 19$ and hence $Sub(G)\geq Sub(G/N)+Sub(N)-1=59+19-1=77$, a contradiction. Thus $n\leq 3$. If $n=3$, then $Sub(N)=2p^2+2p+4$ and only prime satisfying $Sub(N)\leq 18$ is $p=2$. If $n=2$, $Sub(N)=p+3$ and only primes satisfying $Sub(N)\leq 18$ are $p=2,3,5,7,11,13$. Thus we have a very few choice for $|G|=60p^n$. Now, using GAP, one can check that among all perfect groups $G$ of these orders, none of them satisfies $Sub(G)=76$. Thus we must have $n=1$, i.e., $N\cong \mathbb{Z}_p$ and $|G|=60p$.

 If $p=2,3$ or $5$, then $|G|=120,180$ or $300$. Among perfect groups of these orders, only tenable candidate is $G\cong SL(2,5)$ with $Sub(SL(2,5))=76$.

 If $p\geq 7$, then as $N\lhd G$ and $gcd(|N|,|G|/|N|)=1$, by Schur–Zassenhaus theorem, $N$ has a complement in $G$, i.e., $G$ has a subgroup of order $60$. Again, since the subgroups of order $2,3,4$ and $5$ are not normal in $G$, $G$ has at least  $$1+(1+2)+(1+2)+(1+3)+(1+5)+1=18 \mbox{ subgroups},$$ which are not containing $N$. Here the first $1$ denotes the trivial subgroup and the last $1$ denotes the subgroup of order $60$. Thus $G$ has at least $59+18=77$ subgroups, a contradiction. This proves the theorem.\qed 

\section{Supersolvable Enforcing Numbers}\label{ss-enf-numbers}
In this section, we get show that if $22\leq Sub(G)\leq 24$, then $G$ is supersolvable.
\begin{theorem}
    Let $G$ be a non-supersolvable group such that $22\leq Sub(G)\leq 24$, then $|G|$ has exactly $2$ distinct prime factors.
\end{theorem}
\pf Clearly $G$ is not a $p$-group. Also $|G|$ is not square-free. If possible, let $|G|$ has at least three prime factors. \\
{\bf Claim:} $G$ has no normal subgroup of prime order.\\
{\it Proof of Claim:} Let $H$ be a normal subgroup of prime order $r$. Then $G/H$ is non-supersolvable. Then either $Sub(G/H)=10$ and $G/H \cong A_4$ or $Sub(G/H)=15$ and $G/H \cong SL(2,3)$ or $(\mathbb{Z}_2\times \mathbb{Z}_2)\rtimes \mathbb{Z}_9$ or $Sub(G/H)\geq 20$. 

Suppose $Sub(G/H)\geq 20$. If for all prime divisors $p (\neq r)$ of $|G|$, the $p$-subgroups are unique in $G$, then the product $H$ of all Sylow $p$-subgroups ($p\neq r)$ is a cyclic and normal subgroup of $G$. Also $G/H$, being a $r$-group, is supersolvable, which implies $G$ is supersolvable, a contradiction. Thus there exists at least one prime divisor $q\neq r$ of $G$ such that $G$ has at least $1+q$ many $q$-subgroups of $G$ of some particular order $q^k$.  Thus counting the subgroups of $G$ we get $$1+20+(1+q)+1=23+q\geq 25, \mbox{ a contradiction}$$ 
where the first $1$ stands for $\{e\}$, $20$ stands for the subgroup of $G$ containing $H$, $1+q$ gives the number of $q$-subgroups and the last $1$ stands for Sylow subgroup corresponding to the third prime factor.

If $Sub(G/H)=10$ and $G/H \cong A_4$. Then $|G|=2^2\cdot 3\cdot r$ (where $|H|=r\geq 5$). Note that this implies that Sylow-$r$-subgroup $H$ is unique and normal in $G$,  i.e., $n_r=1$.

Therefore, there exists exactly $10$ subgroups of $G$ containing $H$ (including $H$ and $G$). In fact, there are precisely $10$ subgroups of $G$ whose order is divisible by $r$. If the number of Sylow $3$-subgroup, $n_3=1$, we get a cyclic normal subgroup $K$ of order $3r$ and $|G/K|=4$, i.e., $G/K$ is supersolvable, a contradiction. Thus $n_3\geq 4$. If subgroup of order $2$ in $G$ is unique, say $L$, then $|G/L|=6r$, i.e., square-free and hence $G/L$ is supersolvable, a contradiction. Thus we get at least $3$ subgroups of order $2$. Also, as $G$ is solvable, due to existence of Hall subgroups, $G$ has a subgroup of order $12$, say $T_1$. Clearly, $T_1$ is not normal in $G$, as otherwise $G\cong T_1\times H$ and hence $Sub(G)=2\times Sub(T_1)$, a contradiction comparing all groups of order $12$. Let $T_2$ be a conjugate subgroup of $T_1$ of order $12$, i.e., $T_1\cong T_2$.

Clearly $T_1$ can not have any element of order $4$, as otherwise Sylow $2$,$3$,$r$-subgroups are all cyclic and thereby making $G$ supersolvable. Thus $T_1\not\cong \mathbb{Z}_{12},\mathbb{Z}_3\rtimes \mathbb{Z}_4$. If $T_1$ has an element of order $6$, then $G$ has a subgroup $M$ of order $6r$, which implies $G/H\cong A_4$ has a subgroup of order $6$, a contradiction. Thus $T_1\not\cong \mathbb{Z}_6\times \mathbb{Z}_2$. If $T_1\cong D_6$, then $Sub(T_1)=16$ and hence $10+16>24$, a contradiction. Thus $T_1\not\cong D_6$. So, we have $T_1\cong T_2\cong A_4$.

Now, we consider the subgroup $T_1\cap T_2$. Its possible orders are $1,2,3,4,6$. If $|T_1\cap T_2|=6$, then similarly as above, we get a subgroup of order $6$ in $A_4$, a contradiction. If $|T_1\cap T_2|=4$, then $T_1\cup T_2$ contains $8$ subgroups of order $3$. Thus $n_3=1+3s\geq 10$. Hence $$1+10+10+3+1=25>24,$$
where first $1$ stands for the trivial subgroup, first $10$ stands for subgroups of $G$ containing $H$, second $10$ stands for $n_3$, $3$ stands for $2$-order subgroups and the last $1$ stands for $n_2$. If $|T_1\cap T_2|=3$, then $T_1\cup T_2$ contains $7$ subgroups of order $3$ and $6$ subgroups of order $2$. Thus no. of $2$-order subgroups in $G$ is at least $7$. Hence $$1+10+7+7+1=26>24,$$ where first $1$ stands for trivial subgroup, $10$ stands for subgroups of $G/H$, first $7$ stands for $n_3$, the second $7$ stands for $2$-order subgroups and the last $1$ stands for $n_2$. If $|T_1\cap T_2|=2$, then $T_1\cup T_2$ contains $5$ subgroups of order $2$ and $8$ subgroups of order $3$, i.e., $n_3\geq 10$. Hence $$1+10+10+5+1=27>24,$$ where first $1$ stands for the trivial subgroup, first $10$ stands for subgroups of $G/H$, second $10$ stands for $n_3$, $5$ stands for $2$-order subgroups and the last $1$ stands for $n_2$. Similarly, it can be shown that $|T_1\cap T_2|\neq 1$. Thus $Sub(G/H)\neq 10$.

If $Sub(G/H)=15$ and $G/H \cong SL(2,3)$, then $|G|=2^3\cdot 3 \cdot r$. Arguing as in the previous case, we get $n_r=1$, $n_3\geq 4$, at least two Hall subgroups of order $24$ and at least one subgroup each of order $2,2^2$ and $2^3$. Thus we get at least $1+15+4+2+3\geq 25$ subgroups of $G$, a contradiction. 

If $Sub(G/H)=15$ and $G/H\cong (\mathbb{Z}_2\times \mathbb{Z}_2)\rtimes \mathbb{Z}_9$, then $|G|=2^2\cdot 3^2 \cdot r$. Using similar counting argument, we reach a contradiction. Thus the Claim is proved.

From Claim, it follows that the number of subgroups of order $p$, $q$ and $r$ in $G$ are $1+pk_1,1+qk_2$ and $1+rk_3$ respectively, where $k_1,k_2,k_3\geq 1$. As $G$ is solvable, $G$ has at least $3$ Hall subgroups corresponding the set of primes $\{p,q\},\{q,r\}$ and $\{p,r\}$. Again, as $G$ is not nilpotent, $G$ has a maximal subgroup $M$ of prime-power index in $G$ such that $M$ is not normal in $G$. Thus there exists at least two maximal subgroups of order $|M|$, namely $M$ and a conjugate subgroup of $M$. (note that $M$ may be a Hall subgroup) If $M$ is a Hall $\{p,q\}$-subgroup, then the number of Hall $\{p,q\}$-subgroups, $h_m\geq 3$ (by Theorem 9.3.1, \cite{hall-book}). Thus we get at least $5$ subgroups counting the maximal subgroups which are not normal and the Hall subgroups of $G$. Thus 
$$2+1+(1+pk_1)+(1+qk_2)+(1+rk_3)+5=11+pk_1+qk_2+rk_3\leq 24,$$
 where $2$ stands for the trivial subgroup and the group $G$, $1$ stands for a subgroup of prime-squared order (since $|G|$ is not square-free), i.e., 
 \begin{equation}\label{3-prime-factor-eq-1}
     pk_1+qk_2+rk_3\leq 13
 \end{equation}
 From Equation \ref{3-prime-factor-eq-1}, it follows that $|G|$ has exactly $3$ prime factors. Also the only possibilities of $\{p,q,r\}$ are $\{2,3,5\}$ and $\{2,3,7\}$.

We first eliminate the case when $\{p,q,r\}=\{2,3,7\}$. As $2+3+7=12$, i.e., $1$ less than $13$, the only possible options available for $|G|$ is $p^3qr$ or $p^2q^2r$ (otherwise the subgroup count will exceed at least by $2$. Now, by performing an exhaustive search on non-supersolvable groups of orders $p^3qr$ or $p^2q^2r$ where $\{p,q,r\}=\{2,3,7\}$ using GAP, one can see that no such group $G$ exists with $22\leq Sub(G)\leq 24$. 

Now, we deal with the case $\{p,q,r\}=\{2,3,5\}$. Let $|G|=2^\alpha3^\beta 5^\gamma$, where $\alpha+\beta+\gamma\geq 4$ and $\alpha,\beta,\gamma\geq 1$. As $2+3+5=10$, i.e., $3$ less than the upper bound given by Equation \ref{3-prime-factor-eq-1}, we have $\alpha+\beta+\gamma \leq 7$. Moreover, if $\alpha+\beta+\gamma =6$ or $7$, then at least one subgroup of order $p^2$ or $p^3$ is unique and hence normal in $G$. This will give rise to some additional subgroups of order $p^2q$ or $p^3q$, which were not included in the count given by Equation \ref{3-prime-factor-eq-1}. Thus Equation \ref{3-prime-factor-eq-1} will be violated. Hence we have $$4\leq \alpha+\beta+\gamma\leq 5.$$ Now, by performing an exhaustive search on non-supersolvable groups of orders $2^\alpha3^\beta 5^\gamma$ with $4\leq \alpha+\beta+\gamma\leq 5$ using GAP reveals that no such group $G$ exists with $22\leq Sub(G)\leq 24$. Hence the theorem follows.\qed 

\begin{lemma}
    There does not exist any non-supersolvable group $G$ such that $22\leq Sub(G)\leq 24$ and $|G|=p^\alpha q$.
\end{lemma}
\pf Suppose such a group $G$ exists. Clearly $\alpha \geq 2$ and $q\nmid (p-1)$ (otherwise $G$ will be supersolvable). As Sylow $q$-subgroup is cyclic, it is not normal in $G$. Thus the number of Sylow $q$ subgroups in $G$ is $n_q=1+qk\geq p^2$. Again, as Sylow $q$-subgroup is cyclic, the Sylow $p$-subgroup(s) is/are not cyclic. Thus by Theorem A, \cite{min-no-of-sg}, we get at least $(\alpha-1)p+(\alpha+1)$ many $p$-subgroups. Thus counting the number of subgroups of $G$, we get 
\begin{equation}\label{palphaq-eq-1}
    (\alpha-1)p+(\alpha+1)+n_q+1\leq 24
\end{equation}
i.e., $p+p^2\leq 20$, i.e., $p=2$ or $3$. 

If $p=3$, then $n_q=1+qk=9$ implies $q=2$, i.e., $q|(p-1)$, a contradiction. Thus $p=2$. Now, $n_q=1+qk=2^2,2^3$ or $2^4$ (as $n_q\geq 2^5$ will make the total number of subgroups exceed $24$) implies that $q=3,5$ or $7$. 

If $q=3$, then $n_q\geq 4$, then Equation \ref{palphaq-eq-1} yields $\alpha\leq 6$. If $q=5$, then $n_q=16$ and hence from Equation \ref{palphaq-eq-1}, we get $\alpha=2$. If $q=7$, then $n_q=8$ and hence from Equation \ref{palphaq-eq-1}, we get $\alpha \leq 5$.

Now an exhaustive search on non-supersolvable groups of orders $2^\alpha\cdot 3$ with $2\leq \alpha\leq 6$, $2^2\cdot 5$ and $2^\alpha\cdot 7$ with $2\leq \alpha\leq 5$ reveals that no such group $G$ with $22\leq Sub(G)\leq 24$ exists. Hence the theorem holds. 

\begin{theorem}
    Let $G$ be a non-supersolvable group such that $22\leq Sub(G)\leq 24$ and $|G|=p^\alpha q ^\beta$. Then $\alpha+\beta\leq 6$.
\end{theorem}
\pf If possible, let $\alpha+\beta\geq 7$. Clearly both Sylow $p$ subgroup and Sylow $q$ subgroup of $G$ can not be normal in $G$.\\
{\bf Claim 1:} None of the Sylow Subgroups of $G$ are normal.\\
{\it Proof of Claim 1:} Suppose the Sylow $p$-subgroup $S_p \lhd G$ and Sylow $q$-subgroup $S_q \ntriangleleft G$. Moreover $S_p$ is not cyclic (as $S_p$ is cyclic and $G$ is non-supersolvable implies $G/S_p$ is non-supersolvable, a contradiction).

Let $H$ be a $q$-subgroup of $G$ of order $q^k$, where $k \in \{1,2,3,4,5\}$ (if it exists). Then we claim that $H\ntriangleleft G$. Because, if $H\lhd G$, then $G/H$ is non-supersolvable (as $G/H$ and $G/S_p$ are supersolvable implies $G \cong G/(H\cap S_p)$ is supersolvable). Thus $Sub(G/H)\in \{10,15,20,21,22,23,24\}$. If possible, let $Sub(G/H)\geq 20$. Then as $\alpha\geq 2$ and $S_p$ is not cyclic, the number of $p$-subgroups of $G$ is greater than or equal to $p+2$. Thus $1+(p+2)+Sub(G/H)\geq 23+p\geq 25$, a contradiction. Thus $Sub(G/H)=10$ or $15$, i.e., $G/H \cong A_4$ or $SL(2,3)$ or $(\mathbb{Z}_2\times \mathbb{Z}_2)\rtimes \mathbb{Z}_9$, i.e., $|G/H|=2^2\cdot 3$ or $2^3\cdot 3$ or $2^2\cdot 3^2$. Hence $|G|=2^{k+2}\cdot 3$, $2^2\cdot 3^{k+1}$, $2^{k+3}\cdot 3$, $2^3\cdot 3^{k+1}$, $2^{2+k}\cdot 3^2$ or $2^2\cdot 3^{k+2}$, where $k=1,2,3,4,5$. An exhaustive search on non-supersolvable groups $G$ of these orders reveals that $22\leq Sub(G) \leq 24$ do not hold for these groups. Thus $H\ntriangleleft G$, i.e., subgroups of orders $q,q^2,q^3,q^4,q^5$ (if they exist) in $G$ are not unique and we get at least $1+q$ many subgroups of order $q,q^2,q^3,q^4,q^5$ (if they exist)

Since $p^2$ divides $|G|$ and $S_p$ is not cyclic, we have $Sub(S_p)\geq p+3$. If $q^5||G|$, then we should have $$5(1+q)+(p+3)+4+1=5q+p+13\leq 24, \mbox{ i.e., }5q+p\leq 11, \mbox{ a contradiction.} $$
Here subgroups of order $p^\alpha q,p^\alpha q^2,p^\alpha q^3, p^\alpha q^4$ contribute $4$ and $G$ itself contributes $1$ to the sum.

If $q^4||G|$, but $q^5\nmid |G|$, then $\alpha\geq 3$. Since $S_p$ is not cyclic, we have $Sub(S_p)\geq 6$ if $p=2$ and $Sub(S_p)\geq 10$ if $p>2$. Thus, if $p=2$, we get $4(1+q)+6+3+1=4q+14\leq 24$, i.e., $4q \leq 10$, a contradiction. Also, if $p>2$, we get $4(1+q)+10+3+1=4q+18\leq 24$, i.e., $4q\leq 6$, a contradiction.

If $q^3||G|$, but $q^4\nmid |G|$, then $\alpha\geq 4$. Then $Sub(S_p)\geq 11$ if $p=2$ and $Sub(S_p)\geq 14$ if $p>2$. Thus, if $p=2$, we get $3(1+q)+11+2+1=3q+17\leq 24$, i.e., $3q \leq 7$, a contradiction. Also, if $p>2$, we get $3(1+q)+14+2+1=3q+20\leq 24$, i.e., $3q\leq 4$, a contradiction.

If $q^2||G|$, but $q^3\nmid |G|$, then $\alpha\geq 5$. Then $Sub(S_p)\geq 14$ if $p=2$ and $Sub(S_p)\geq 18$ if $p>2$. Thus, if $p=2$, we get $2(1+q)+14+1+1=2q+18\leq 24$, i.e., $2q \leq 6$, i.e., $q=3$. Note that this holds if $\alpha=5$ and hence $|G|=2^5\cdot 3^2$. An exhaustive search on non-supersolvable groups $G$ of order $32\cdot 9$ reveals that $Sub(G) \geq 56$, a contradiction. Also, if $p>2$, we get $2(1+q)+18+1+1=2q+22\leq 24$, i.e., $2q\leq 2$, a contradiction.

Hence Claim 1 holds, i.e., Sylow subgroups of $G$ are not normal. Thus, we have $p+Sub(S_p)+q+Sub(S_q)-1+1\leq 24$, i.e., 
\begin{equation}\label{a+bleq6-equation}
    Sub(S_p)+Sub(S_q)+(p+q)\leq 24.
\end{equation}

We also note that both the Sylow subgroups of $G$ can not be cyclic, as that would imply that $G$ is supersolvable. On the other hand, in the next claim, we show that both of them can not be non-cyclic.\\
{\bf Claim 2:} Both the Sylow Subgroups of $G$ can not be non-cyclic.\\
{\it Proof of Claim 2:} Since $\alpha+\beta \geq 7$ and $\alpha,\beta \geq 2$, we have the following cases:

If $\alpha=2$, then $\beta\geq 5$ and hence from Equation \ref{a+bleq6-equation}, we get $(p+3)+(4q+6)+(p+q)\leq 24$, i.e., $2p+5q\leq 15$, a contradiction.

If $\alpha=3$, then $\beta\geq 4$ and Equation \ref{a+bleq6-equation} yields
$$\left \lbrace \begin{array}{cl}
    \mbox{for }p=2, & 6+(3q+5)+(2+q)=4q+13\leq 24, \mbox{ i.e., }4q\leq 11 \\
    \mbox{for }p>2, & (2p+4)+(3q+5)+(p+q)=3p+4q+9 \leq 24, \mbox{ i.e., }3p+4q\leq 15, 
\end{array}\right.$$
none of which can hold.

The cases $\alpha=4,\beta\geq 3$ and $\alpha=5,\beta \geq 2$ also give rise to contradiction, as in the above two cases. Hence Claim 2 holds.

So, without loss of generality, we assume that $S_p$ is cyclic and $S_q$ is not cyclic.

If any subgroup $K$ of order $p^k$ with $1\leq k \leq 4\leq \alpha-1$ (provided it exists) is normal and as $H$ is cyclic, $G/H$ is non-supersolvable.  Therefore $Sub(G/H) \in \{10,15,20,21,22,23,24\}$. Now proceeding as above, we get a contradiction.  Thus subgroups of order $p^k$ with $1\leq k \leq 4$  (if it exists) are not normal, and hence not unique in $G$. Thus we have at least $1+p$ subgroups for each power of $p$ dividing $|G|$. 

Thus we can replace $Sub(S_p)+p$ in Equation \ref{a+bleq6-equation} by $\alpha(1+p)+1$ to get 
\begin{equation}\label{a+bleq6-equation-2}
    \alpha(1+p)+1+Sub(S_q)+q\leq 24
\end{equation}

If $\alpha\geq 5$ and as $\beta\geq 2$, we have $Sub(S_q)\geq q+3$. Thus from Equation \ref{a+bleq6-equation-2}, we get $5(1+p)+1+(q+3)+q\leq 24$, i.e., $5p+2q\leq 15$, a contradiction.

If $\alpha=4$, then $\beta\geq 3$, i.e., $Sub(S_q)\geq 6$ if $q=2$ and $Sub(S_q)\geq 2q+4$, if $q>2$. Thus from Equation \ref{a+bleq6-equation-2}, if $q=2$, we have $4(1+p)+1+6+2=4p+13\leq 24$, i.e., $4p\leq 11$, a contradiction.
Similarly, if $q>2$, we get $4(1
+4)+q=4p+3q+9\leq 24$, i.e., $4p+3q\leq 15$, a contradiction.

If $\alpha=3$, then $\beta\geq 4$, i.e., $Sub(S_q)\geq 3q+5$. Thus from Equation \ref{a+bleq6-equation-2}, we have $3(1+p)+1+(3q+5)+q=3p+4q+9\leq 24$, i.e., $3p+4q\leq 15$, a contradiction.

If $\alpha=2$, then $\beta\geq 5$, i.e., $Sub(S_q)\geq 4q+1$. Thus from Equation \ref{a+bleq6-equation-2}, we have $2(1+p)+1+(4q+1)+q=2p+5q+9\leq 24$, i.e., $2p+5q\leq 15$, a contradiction.

Thus $\alpha+\beta \geq 7$ can not hold and the theorem follows.

\begin{lemma}
    There does not exist any non-supersolvable group $G$ such that $22\leq Sub(G)\leq 24$ and $|G|=p^\alpha q^2$ for $\alpha=2,3$ or $4$.
\end{lemma}
\pf Suppose such a group $G$ exists. \\
{\bf Case 1:} (Sylow $q$-subgroup $S_q$ is cyclic).  So $S_q$ is not normal in $G$. Thus the number of Sylow $q$-subgroup of $G$ is $n_q=1+qk\geq p$ and $n_q|p^\alpha$. Again, as $S_q$ is cyclic, the Sylow $p$-subgroup $S_p$ is not cyclic. Thus counting subgroups of $G$, we get 
\begin{equation}\label{palphaq2-eq-1}
    Sub(S_p)+1+n_q+1\leq 24
\end{equation}

For $\alpha=4$, we have $Sub(S_p)\geq 3p+5$. Thus from Equation \ref{palphaq2-eq-1}, we get $(3p+5)+1+n_q+1\leq 24$. As $n_q\geq p$, we get $4p\leq 17$, i.e., $p=2$ or $3$. If $p=2$, then $n_q=1+qk|16$ which implies $q=3$ or $5$. If $q=5$, then $n_q=1+5k=16$, which exceeds the count of subgroups. Thus the only possible order of $G$ is $2^4\cdot 3^2$.

For $\alpha=3$, using Equation \ref{palphaq2-eq-1} and arguing as above, we get the only possible orders of $G$ as $3^3\cdot 2^2, 2^3\cdot 3^2, 2^3\cdot 7^2$ and $5^3\cdot 2^2$.

Similarly for $\alpha=2$, the possible orders are $2^2\cdot 3^2,2^2\cdot 5^2$ and $2^2\cdot 7^2$.

Now an exhaustive search on non-supersolvable groups of above orders reveals that no such group $G$ with $22\leq Sub(G)\leq 24$ exists.

{\bf Case 2:} (Sylow $q$-subgroup $S_q$ is not cyclic). As  Sylow $p$-subgroup $S_p$ is not simultaneously cyclic and normal in $G$, similar counting arguments leaves only finitely many possible orders of $G$ which can be computationally checked not have $22\leq Sub(G)\leq 24$.\qed

\begin{lemma}\label{p3q3}
    There does not exist any non-supersolvable group $G$ of order $p^3q^3$ such that $22\leq Sub(G)\leq 24$.
\end{lemma}
\pf Suppose such a group $G$ exists. As at least one of the Sylow subgroups is not cyclic (say $S_p$), and $S_q$ is either not cyclic or not normal in $G$, we have 
\begin{equation}
    Sub(S_p)+(1+1+1+q)+2\leq 24
\end{equation}
If $p=2$, this gives $6+(3+q)+2=q+11\leq 24$, i.e., $q=2,3,5,7,11$ or $13$. If $p>2$, then the above equation gives $(2p+4)+(3+q)+2=2p+q\leq 15$. In both cases, we have a few choices for $|G|$ and can be exhaustively checked not to have $Sub(G)$ lying between $22$ and $24$. \qed

Combining all the lemmas and theorems of this section, we get the following theorem:

\begin{theorem}\label{22-24-SS-Theorem}
    Let $G$ be a group such that $22\leq Sub(G)\leq 24$, then $G$ is supersolvable.
\end{theorem}

\begin{corollary}
    Let $G$ be a nilpotent group such that $Sub(G)=22$, then $G$ is one of the following form:
    \begin{itemize}
    \item $G$ is cyclic and $G\cong \mathbb{Z}_{p^{21}}$ or $G\cong \mathbb{Z}_{pq^{10}}$, where $p,q$ are distinct primes.
    \item $G$ is a non-cyclic $p$-group and $G$ is isomorphic to one of the following groups:
    \begin{itemize}
        \item GAP ID $(32,3)$: $\mathbb{Z}_8 \times \mathbb{Z}_4$.
        \item GAP ID $(32,4)$: $\mathbb{Z}_8 \rtimes \mathbb{Z}_4$.
        \item GAP ID $(32,12)$: $\mathbb{Z}_4 \rtimes \mathbb{Z}_8$.
        \item GAP ID $(361,2)$: $\mathbb{Z}_{19} \times \mathbb{Z}_{19}$.
        \item GAP ID $(729,93)$: $\mathbb{Z}_{243} \times \mathbb{Z}_3$.
        \item GAP ID $(729,94)$: $\mathbb{Z}_{243} \rtimes \mathbb{Z}_3$.
    \end{itemize}
    \item $G$ is non-cyclic non-$p$-group and $G$ is isomorphic to one of the following groups, where $p$ is an odd prime:
    \begin{itemize}
        \item $\mathbb{Z}_p \times \mathbb{Z}_8\times \mathbb{Z}_2$.
        \item $\mathbb{Z}_p \times Q_{16}$.
        \item $\mathbb{Z}_p \times M_{16}$.
    \end{itemize}
    \end{itemize}
\end{corollary}
\pf The proof follows from the fact that $Sub(G)=22=2\cdot 11$ and $G$ is nilpotent implies that $|G|$ has at most two prime factors and all Sylow subgroups of $G$ are normal in $G$. We omit the details for brevity.

\section{Nilpotent Enforcing Numbers}\label{nil-enf-numbers}
In this section, we show that if $Sub(G)=23$, then $G$ is nilpotent and $G$ is isomorphic to one of the seven groups given in Theorem \ref{Sub(G)=23-main-theorem}.
\begin{lemma}\label{Sub=23-numberofprimes}
    Let $G$ be a non-nilpotent group such that $Sub(G)=23$. Then $|G|=pqrs$ or $|G|$ has at most three distinct prime factors, where $p,q,r,s$ are distinct primes.
\end{lemma}
\pf Since $Sub(G)=23$, by Theorem \ref{22-24-SS-Theorem}, $G$ is supersolvable and hence a CLT group. Also, since $G$ is non-nilpotent, $G$ is not a $p$-group and $G$ has at least one Sylow subgroup which is not normal. Let $|G|=p^{\alpha_1}_1p^{\alpha_2}_2\cdots p^{\alpha_k}_k$
and Sylow $p_1$-subgroup be not normal in $G$. Then we have $(\alpha_1+1)(\alpha_2+1)\cdots (\alpha_k+1)+p_1\leq 23$, i.e., 
\begin{equation}\label{23-eq-1}
  (\alpha_1+1)(\alpha_2+1)\cdots (\alpha_k+1)\leq 21, \mbox{ where }\alpha_i\geq 1.  
\end{equation}
As all integers less than $22$, except $16$, can not be expressed as product of more than $3$ factors $>1$, either $G$ has at most $3$ distinct prime factors or if $(\alpha_1+1)(\alpha_2+1)\cdots (\alpha_k+1)=16$, then $|G|=pqrs$. \qed

\begin{lemma}\label{Sub=23-2primes}
    There does not exist any non-nilpotent group $G$ such that $Sub(G)=23$ and $|G|=p^\alpha q^\beta$. 
\end{lemma}
\pf Let $G$ be a non-nilpotent group such that $Sub(G)=23$ and $|G|=p^\alpha q^\beta$. If $\alpha+\beta\geq 11$ with $\alpha,\beta\geq 1$, then the minimum value of $(\alpha+1)(\beta+1)\geq 22$ which contradicts Equation \ref{23-eq-1}. Thus $\alpha+\beta\leq 10$. If $7\leq \alpha+\beta\leq 10$, the only possible options obeying Equation \ref{23-eq-1} is $|G|=p^9q,p^8q,p^7q,p^6q,p^6q^2$ and $p^5q^2$. We first show that such orders of $G$ is not possible.

Suppose $|G|=p^\alpha q$, where $6\leq \alpha\leq 9$. If $q>p$, then the Sylow $q$-subgroup is normal in $G$ and Sylow $p$-subgroups are not normal in $G$, i.e., $1<n_p=1+pk|q$, i.e., $n_p=q$. If $S_p$ is not cyclic, then $Sub(S_p)\geq (\alpha-1)p+(\alpha+1)\geq 5p+7$. Thus we have $5p+7+(q-1)+(\alpha+1)\leq 23$, i.e., $5p+q+13\leq 23$, i.e., $5p+q\leq 10$, a contradiction. Thus $S_p$ must be cyclic and we have $2(1+\alpha)+(q-1)\leq 23$, i.e., $q\leq 10$. So the only choice for $q$ is $3,5,7$ and $p<q$. Similarly, if $p>q$, then the Sylow $p$-subgroup is normal in $G$ and Sylow $q$-subgroups are not normal in $G$, i.e., $1<n_q=1+qk|p^\alpha$, i.e., $n_q\geq q$. In the same, way it can be shown that $p$ is either $3$ or $5$ or $7$. Thus, in any case, $p,q \in \{2,3,5,7\}$ and $6\leq \alpha\leq 9$. So, we need to look only among non-nilpotent, supersolvable groups $G$ with $|G|=p^\alpha q$ and both Sylow subgroups cyclic.  It can be checked using GAP that no such group has $Sub(G)=23$. Hence $|G|\neq p^9q,p^8q,p^7q,p^6q$.

Suppose $|G|=p^\alpha q^2$, where $5\leq \alpha\leq 6$. Using similar counting techniques as above, one can show that $p,q\in \{2,3,5\}$ and then the rest can be checked with GAP. Thus $|G|\neq p^5q^2,p^6q^2$.

So, we have $\alpha+\beta\leq 6$. If $|G|=p^\alpha q^\beta$ with $3\leq \alpha+\beta\leq 6$. Without loss of generality, let $p>q$. Then $S_p \lhd G$ and $S_q$ is not normal in $G$ and $n_q\geq p$. Thus we have $(1+\alpha)(1+\beta)+(p-1)\leq 23$, i.e., $\alpha\beta+(\alpha+\beta)+p\leq 23$. So, we have a very few choices for $p,q$ and with small exponents $\alpha,\beta$. Using GAP, one can check that such group $G$ does not exist.

So, lastly we are left with the case when $|G|=pq$ where $G\cong \mathbb{Z}_p\rtimes \mathbb{Z}_q$ where $p>q$. In this case, $Sub(G)=p+3$ which implies $p+3=23$, i.e., $p=20$, which is not a prime.

Hence the lemma follows.\qed
 
\begin{lemma}\label{Sub=23-3primes}
    There does not exist any non-nilpotent group $G$ such that $Sub(G)=23$ and $|G|=p^\alpha q^\beta r^\gamma$. 
\end{lemma}
\pf Let $G$ be a non-nilpotent group such that $Sub(G)=23$ and $|G|=p^\alpha q^\beta r^\gamma$. If $\alpha+\beta+\gamma\geq 7$, then the minimum value of $(\alpha+1)(\beta+1)(\gamma+1)$ is $24$, which exceeds the upper bound given by Equation \ref{23-eq-1}. So, we must have $\alpha+\beta+\gamma\leq 6$. 

If $\alpha+\beta+\gamma=6$, then the only possible order of $G$ satisfying Equation \ref{23-eq-1} is $p^4qr$. Similarly, if $\alpha+\beta+\gamma=5$, the only possible orders of $G$ are $p^3qr$ and $p^2q^2r$. For $3\leq \alpha+\beta+\gamma\leq 4$, we have $|G|=pqr$ or $p^2qr$.

If $|G|=p^4qr$, as $G$ is Lagrangian, we get at least $5\cdot 2\cdot 2=20$ subgroups, counting subgroups of each order exactly once. Also note that two of the Sylow subgroups can not be non-normal in $G$, as it would exceed the count. So, exactly one of the Sylow subgroup is not normal in $G$ and other two are normal in $G$. Let $H,K$ be normal Sylow subgroups of $G$ and $L$ be a non-normal Sylow subgroup of $G$. Then $L$ has at least $2$ conjugate subgroups which are not included in the count of $20$. Moreover $HL$ is a Hall subgroup of $G$. If it is a unique Hall subgroup of order $HL$, then $G=HL\times K$ and $23=Sub(G)=Sub(HL)\cdot Sub(K)$, a contradiction. So there are at least $2$ more Hall subgroups of order $|HL|$ which are not included in the count of $20$. So, in total we get at least $4=2+2$ subgroups other than the count of $20$, i.e., $Sub(G)\geq 24$, a contradiction.

If $|G|=p^2q^2r$, as $G$ is Lagrangian, we get at least $3\cdot 3\cdot 2=18$ subgroups. If exactly two of the Sylow subgroups of $G$ are non-normal in $G$, then the corresponding primes must be $2$ and $3$ and the total count becomes $23$, i.e., all subgroups apart from the Sylow $2$-subgroups and Sylow $3$-subgroups are of unique in $G$. Let $H$ be the unique Hall $\{2,3\}$-subgroup of $G$ and $K$ be the normal Sylow subgroup corresponding to the third prime. Thus $G\cong H \times K$ and we get a contradiction, as above. So, exactly one of the Sylow subgroups is not normal in $G$ and other two are normal in $G$. Let $H,K$ be normal Sylow subgroups of $G$ and $L$ be a non-normal Sylow subgroup of $G$. Then $L$ has at least $2$ conjugate subgroups which are not included in the count of $18$. Moreover $HL$ and $KL$ are Hall subgroup of $G$.
If any of them is unique, then $G$ can be expressed as direct product of two subgroups of coprime order and we get a contradiction as above. So, there are at least $4$ more Hall subgroups ($2$ each of orders $|HL|$ and $|KL|$) which are not included in the count of $18$. So, in total we get at least $6=2+4$ subgroups other than the count of $18$, i.e., $Sub(G)\geq 24$, a contradiction.

If $|G|=pqr$ with $p<q<r$, then $n_r=1$, i.e., the Sylow $r$-subgroup, $R$ and the Hall subgroup $H_{q,r}$ of order $qr$ are normal in $G$, i.e., $G\cong H_{q,r}\rtimes \mathbb{Z}_p$ and $n_p \in \{q,r,qr\}$, i.e., $n_p\geq q$. If $H_{q,r}\cong \mathbb{Z}_r \rtimes \mathbb{Z}_q$, then $Sub(H_{q,r})=r+3$. Now consider a Hall $\{p,q\}$-subgroup, $H_{p,q}$. If it is normal in $G$, then $G \cong H_{p,q}\times R$. This contradicts that $Sub(G)=23$ is prime. Thus $G$ has at least $3$ Hall $\{p,q\}$-subgroups and they are not normal in $G$. Now, counting the subgroups of $G$, we get (the last two $1$'s are for Hall $\{p,r\}$-subgroup and $G$ itself) $$23=Sub(G)\geq Sub(H)+n_p+3+1+1=(r+3)+q+3+1+1, \mbox{ i.e., }q+r\leq 15.$$
This leaves very few choices for $(p,q,r)$ and can be exhaustively checked in GAP, i.e., no such groups have exactly $23$ subgroups. Thus $H_{q,r}\cong \mathbb{Z}_{qr}$ and $Sub(H_{q,r})=4$. If $n_p=r$ or $qr$, similarly we get $r$ or $qr\leq 14$. Again, using GAP, one can check that no such group exists. So, we must have $n_p=q$. Now, we try to evaluate $n_q$. Let $T$ be the unique subgroup of order $q$ in $H_{q,r}$. If $T'$ is any other subgroup of order $q$ in $G$, then $|T'\cap H_{q,r}|=1$ and $T'H_{q,r}=G$ (as $H_{q,r}$ is a normal, maximal subgroup of $G$). Thus $|G|=|T'H_{qr}=q^2r>pqr$, a contradiction. Thus $T$ is the unique subgroup of order $q$ in $G$, i.e., $n_q=1$. Now, we consider a Hall $\{p,r\}$-subgroup $H_{p,r}$ of $G$. It is either cyclic or isomorphic to $\mathbb{Z}_r \rtimes \mathbb{Z}_p$. In the later case, $H_{p,r}$ contains $r>q$ subgroups of order $p$, a contradiction. Thus $H_{p,r}$ is cyclic. Similarly, if $H_{p,q}$ is a Hall $\{p,q\}$-subgroup of $G$ and $H_{p,q} \cong \mathbb{Z}_q\rtimes \mathbb{Z}_p$, then $H_{p,q}$ contains all the $q$ many subgroups of order $p$ in $G$. Let $P$ be the subgroup generated by all elements of order $p$ in $G$. Clearly $P$ is a characteristic subgroup of $G$ and $\mathbb{Z}_p \subsetneq P \subseteq H_{p,q}$, i.e., $H_{p,q}=P\lhd G$. However we have proved earlier that Hall $\{p,q\}$-subgroups are not normal in $G$. Thus $H_{p,q}$ is also a cyclic group. Hence all proper subgroups of $G$ are cyclic, i.e., $G$ is a finite minimal non-cyclic groups. Such groups are classified in \cite{miller-moreno} and all of their orders have at most two distinct prime factors, a contradiction.

Similarly, it can be shown that $|G|=p^2qr$ and $p^3qr$ leads to contradiction. We omit the details of calculations for brevity. Hence the lemma holds.\qed


\begin{lemma}\label{Sub=23-pqrs}
    There does not exist any non-nilpotent group $G$ such that $Sub(G)=23$ and $|G|=pqrs$. 
\end{lemma}  
\pf Without loss of generality, let $s$ be the largest prime factor of $|G|$. As $G$ is supersolvable and non-nilpotent, the Sylow $s$-subgroup is normal in $G$ with $n_s=1$ and at least one Sylow subgroup is not normal in $G$. Moreover, as $G$ is Lagrangian, we get at least $16$ subgroups counting the subgroups of each possible order exactly once. If $n_p,n_q,n_r$ are all greater than $1$, then we get $p+q+r\geq 8$ subgroups of $G$ which are not counted in $16$ and hence $Sub(G)\geq 16+8=24$, a contradiction. So, we assume that $n_r=1$. Now, two cases may arise.

{\bf Case 1:} $n_p,n_q\neq 1$ and $n_r=n_s=1$. As $n_p+n_q\geq 2+p+q$, we get at least $p+q$ subgroups of $G$ apart from the count of $16$. Hence we must have $p+q\leq 7$, which implies $\{p,q\}=\{2,3\}$ or $\{2,5\}$. If $\{p,q\}=\{2,3\}$ and $n_3=4$, then $4|2rs$ which is square-free, a contradiction. Thus $n_3\geq 7$ and hence $n_2+n_3\geq 2+2+6$, i.e., we get at least $2+6=8$ subgroups apart from the count of $16$ subgroups, a contradiction. Thus $\{p,q\}=\{2,5\}$. As $n_5>1$, the only option is $n_5=6$ and $n_2=3$. Also as $n_2,n_5$ divides $|G|$, we must have $|G|=2\cdot 3 \cdot 5\cdot s=30s$. 

So, we already have $16$ (number of one subgroups of each possible order) + $5$ (extra Sylow $5$-subgroups of order $5$) + $2$ (extra Sylow $2$-subgroups of order $2$)=$23$ subgroups of $G$. So subgroups of all other orders, except $2$ and $5$ are unique in $G$. Now, consider the Hall subgroups $K$ and $L$, subgroups of order $10$ and $3s$ respectively in $G$. As they are unique, $K,L\lhd G$ and are of co-prime orders in $G$ with trivial intersection, we have $G\cong K\times L$ and $23=Sub(G)=Sub(K)\cdot Sub(L)$, a contradiction. Hence Case 1 can not occur.

{\bf Case 2:} $n_p\neq 1$ and $n_q=n_r=n_s=1$. As we already have at least $16$ subgroups, $p=2,3,5$ or $7$. We rule out each of these cases separately. If $p=7$, then $n_7=8$ and the count becomes exactly $23$, i.e., subgroups of all other orders, except $7$ are unique in $G$ and using arguments as above $G$ can be shown to be direct product of two cyclic groups of co-prime order, a contradiction. If $p=5$, then $n_5=6$ and we get at least $21$ subgroups. Now, out of the Hall subgroups of $G$ of orders $5q,5r,5s,5qr,5qs,qr,qs,rs,qrs$ at most one is non-unique and hence $G$ can be expressed as direct product of two non-trivial subgroups $H$ and $K$ of $G$ with co-prime order. Hence $23=Sub(G)=Sub(H)\times Sub(K)$, a contradiction. If $p=3$, then $n_3=4$ or $7$. However as $4\nmid qrs$, we have $n_3=7$ and we get we get at least $22$ subgroups. So we need one more subgroup to get $23$. However as any subgroup is a Hall subgroup in $G$ and the number of Hall subgroups of any order is either $1$ or $\geq 3$, the count exceeds $23$, a contradiction. If $p=2$, then $n_2=3,5$ or $7$, and $|G|=2qrs$. As $n_2|qrs$, without loss of generality, let $n_2=q$. Let $K$ be a subgroup of order $2r$ and $L$ be a subgroup of order $qs$ in $G$. Clearly $L\lhd G$. If $K$ is the unique subgroup of order $2r$, then $G \cong K\times L$, a contradiction as above. So, as $K$ is a Hall subgroup, $G$ must have at least $3$ Hall subgroups of order $2r$. Arguing similarly, $G$ must have at least $3$ Hall subgroups, each of orders $2s$ and $2rs$. Thus, the total count of subgroups is $\geq 16+2(\mbox{Extra Sylow }2-\mbox{subgroups})+2(\mbox{Extra Hall }2r-\mbox{subgroups})+2(\mbox{Extra Hall }2s-\mbox{subgroups})+2(\mbox{Extra Hall }2rs-\mbox{subgroups})=24$, a contradiction. Hence Case 2 can not occur. Hence the lemma holds.\qed 

From Lemmas \ref{Sub=23-numberofprimes}, \ref{Sub=23-2primes}, \ref{Sub=23-3primes} and \ref{Sub=23-pqrs}, we can say that if $G$ is a group with $Sub(G)=23$, then $G$ is nilpotent. In the next theorem, we fully classify such groups.

\begin{theorem}\label{Sub(G)=23-main-theorem}
    Let $G$ be a non-cyclic group such that $Sub(G)=23$, then $|G|=16,81$ or $256$ and $G$ is isomorphic to one of the $7$ groups given by:
    \begin{enumerate}
        \item GAP ID $(16,3)$: $(\mathbb{Z}_4 \times \mathbb{Z}_2)\rtimes \mathbb{Z}_2$.
        \item GAP ID $(16,13)$: $(\mathbb{Z}_4 \times \mathbb{Z}_2)\rtimes \mathbb{Z}_2$.
        \item GAP ID $(81,2)$: $\mathbb{Z}_9 \times \mathbb{Z}_9$.
        \item GAP ID $(81,4)$: $\mathbb{Z}_9 \rtimes \mathbb{Z}_9$.
        \item GAP ID $(81,10)$: $(\mathbb{Z}_3 \times \mathbb{Z}_3)\cdot (\mathbb{Z}_3 \times \mathbb{Z}_3)$.
        \item GAP ID $(256,537)$: $\mathbb{Z}_{128}\times \mathbb{Z}_2$.
        \item GAP ID $(256,538)$: $\mathbb{Z}_{128}\rtimes \mathbb{Z}_2$.
    \end{enumerate}
\end{theorem}
\pf As $Sub(G)=23$, $G$ is nilpotent and hence $G$ is isomorphic to the direct product of its Sylow subgroups $S_{p_1},S_{p_2},\ldots,S_{p_k}$(say). Thus $$23=Sub(G)=Sub(S_{p_1})\cdot Sub(S_{p_2})\cdots Sub(S_{p_k}).$$ As $23$ is prime, it follows that $G$ is a $p$-group, say $|G|=p^\alpha$. As $G$ is non-cyclic, we have $23=Sub(G)\geq (\alpha-1)p+(\alpha+1)$.

If $\alpha\geq 9$, then $(\alpha-1)p+(\alpha+1)\geq 26$, a contradiction. Hence we have $\alpha\leq 8$. If $\alpha=7$ or $8$, we get $p=2$. If $\alpha=5$ or $6$, then possible values of $p$ are $2$ and $3$. If $\alpha=4$, $p$ is one of $2,3$ or $5$. If $\alpha=3$, then $p \in \{2,3,5,7\}$. If $\alpha=2$, then $Sub(G)=p+3=23$, i.e., $p=20$, which is not a prime. Thus we need the look for groups of order $2^\alpha$ with $3\leq \alpha \leq 8$, $3^\alpha$ with $3\leq \alpha \leq 6$, $5^\alpha$ with $3\leq \alpha \leq 4$ and $7^3$. 

On checking all non-cyclic groups of these orders (except $2^8=256$) using GAP, it reveals that $G$ is isomorphic to one of the first five groups given in the statement of the theorem. 

Now, we deal with the case when $|G|=256$. As there are more than $56000$ non-isomorphic groups of order $256$, to overcome the computational hindrance, we narrow down the search to a lesser number of groups. We start by noting a computational result on non-cyclic groups of order $128$ obtained using GAP: if $G$ is a non-cyclic group of order $128$, then $Sub(G)\geq 36$ except when $G\cong \mathbb{Z}_{64}\times \mathbb{Z}_2$ or $\mathbb{Z}_{64}\rtimes \mathbb{Z}_2$ with $Sub(G)=20$.

Let $G$ be a non-cyclic group of order $256$ with $Sub(G)=23$. Then there exist a normal non-cyclic subgroup $H$ of $G$ of order $2^7=128$. If $H$ is not isomorphic to $\mathbb{Z}_{64}\times \mathbb{Z}_2$ or $\mathbb{Z}_{64}\rtimes \mathbb{Z}_2$, then we have $Sub(G)> Sub(H)\geq 36$, a contradiction. So $H$ must be isomorphic to $\mathbb{Z}_{64}\times \mathbb{Z}_2$ or $\mathbb{Z}_{64}\rtimes \mathbb{Z}_2$ and $Sub(H)=20$. Note that $H$ is a maximal subgroup of $G$. Moreover, $H$ is not the unique maximal subgroup of $G$, as otherwise $G$ would be cyclic. Thus $G$ has $1+2k$ subgroups of order $128$ with $k\geq 1$, out of which one is $H$. Hence $23=Sub(G)\geq Sub(H)+1+2k=21+2k$ and thus $k=1$. So, $G$ has exactly three subgroups of order $128$ out of which one is $H$. Moreover $G$ has rank either $2$ or $3$. Now there are $6730$ groups of order $256$ with rank $2$ or $3$, namely with GAP ID: $2$ to $6731$. We perform our search for $G$ on these groups with a further restriction that it has exactly three maximal subgroups and $Sub(G)=23$. This modified search is run in GAP to get exactly two candidates, namely the last two groups mentioned in the statement of the theorem. \qed


\section{An Application to Comaximal Subgroup Graph of a Group}
The comaximal subgroup graph $\Gamma(G)$ of a group $G$, introduced in \cite{akbari}, is a graph whose set of vertices are non-trivial proper subgroups of $G$ and two vertices $H$ and $K$ are adjacent if $HK=G$. For more details regarding this graph, one can refer to \cite{comaximal-1}, \cite{comaximal-2}. Note that the order of $\Gamma(G)$, i.e., the number of vertices of $\Gamma(G)$ is $Sub(G)-2$. As with many other graphs defined on groups, the most interesting and challenging question is to draw inference about the group from its graph. For example, in \cite{comaximal-2}, it is shown that if the independence number of $\Gamma(G)$ is less than $9$, then $G$ is solvable. On the same note, using our results, one can readily check the solvability, supersolvability, nilpotency of the group just by knowing the order of the graph. For example, using Theorem \ref{Sub<59-solvable}, one can say that if the order of $\Gamma(G)$ is less than $57$, then $G$ is solvable.

\section{Conclusion and Open Issues}
In this paper, we studied the solvability of a group $G$ based on its number of subgroups $Sub(G)$. We also demonstrated an application of it in context of comaximal subgroup graph. However, there are a few interesting observations which we believe to be true but could not prove yet. We conclude by listing some of these observations and directions for further research.
\begin{enumerate}
        \item Apart from the speculations marked in box shown in Figure \ref{fig:enforcing_numbers}, some other observations which were made using GAP for groups of order upto $360$ are as follows:
        \begin{itemize}
            \item If $Sub(G)=26,35,37,52,53,55,58,69$ or $73$, then $G$ is Lagrangian, i.e., converse of Lagrange's theorem holds in $G$.
            \item If $G$ is solvable and $Sub(G)=59$, then $G$ is nilpotent.
        \end{itemize}
	\item As mentioned in previous section, knowing the order of $\Gamma(G)$ may help us to get some information about the group $G$. However, for this we are using only the order of the graph and ignoring other graph parameters. For example, if order of $\Gamma(G)=Sub(G)-2=20$, then by Theorem \ref{22-24-SS-Theorem}, $G$ is supersolvable. But can we say something more about the group looking at the adjacencies of the graph.
\end{enumerate}

\section*{Acknowledgement}
The authors acknowledge the funding of DST-FIST grant $SR/FST/MS-I/2019/41$, Govt. of India. In addition, the first author acknowledge the DST-SERB-MATRICS grant (Sanction no. $MTR/2022/000020$). The authors are also thankful to Roudra Ghosal and Megha Das, as the topic of classifying nature of groups from its number of subgroups was initiated in their Master's thesis and was the starting point of the present work.
	

\end{document}